\documentclass[12pt]{article}

\usepackage{amssymb}

\usepackage[utf8]{inputenc}
\usepackage[T1]{fontenc}
\usepackage{pgf,pgfarrows,pgfnodes,pgfautomata,pgfheaps}
\usepackage{amscd}
\usepackage{amsfonts}
\usepackage{latexsym}
\usepackage{graphicx}
\usepackage{amsthm}
\usepackage{amsmath}
\usepackage{placeins}
\usepackage{float}
\usepackage{pdflscape}

\def\det{\mathrm{det}}

\theoremstyle{definition}

\newtheorem{definicja}{Definicja}[section]

\newtheorem{theorem}[definicja]{Theorem}

\newtheorem{corollary}[definicja]{Corollary}

\newtheorem{lemma}[definicja]{Lemma}
\newtheorem{remark}[definicja]{Remark}
\newtheorem{example}[definicja]{Example}
\def \det{\operatorname{det}}

\begin{document}

\begin{center}
\Large
{\bf A note on algorithmic approach to inverting multivariate formal
power series}\\
\vspace{20pt}
\normalsize
EL\.ZBIETA ADAMUS \\ Faculty of Applied Mathematics, \\ AGH University of Krakow \\
al. Mickiewicza 30, 30-059 Krak\'ow, Poland \\
e-mail: esowa@agh.edu.pl \\
\end{center}

\begin{abstract}
In \cite{ABCH} an effective algorithm for inverting polynomial automorphisms in several variables was proposed. 
We extend its application to the case of multivariate formal power series over a field of arbitrary characteristic and
illustrate the proposed approach with some examples. We establish several properties of the algorithm when applied to mappings with polynomial inverse.  
\end{abstract}

\section{Introduction}

Let $K$ denote a field of an arbitraty characteristic and let $F=(F_1, \ldots, F_n) :K^n \rightarrow K^n$ be a polynomial mapping, i.e. $F_i \in K[X]$, where $X=(X_1, \ldots, X_n)$ for $i=1, \ldots, n$.  $F$ is a \emph{polynomial automorphism} if it  is invertible ant its inverse is also a polynomial mapping.

Assume that the jacobian of $F$ is a nonzero constant.  Such polynomial mapping is called a \emph{Keller map}.
By a linear change one may assume that $F$ is of the form $F(X)=X+H(X)$, where $X=(X_1, \ldots, X_n)$ and $H=(H_1, \ldots, H_n)$ such that for each $H_i$ the order of vanishing is greater than $1$ and the jacobian $\det(J_F)=1$. 
For mapping $F(X)=X+H(X)$ with $H(X)$ homogeneous the jacobian condition is equivalent to jacobian matrix $J_H$ being nilpotent. 

The famous Jacobian Conjecture states that if $K$ is a field of characteristic zero, then every Keller map
$F:K^n \rightarrow K^n$ is a polynomial automorphism. To resolve the conjecutre it is enough to consider mappings of the form
$F(X)=X+H(X)$, where each $H_i$ is homogeneous of degree $3$.

Let us consider mappings of the form
  \begin{equation}
   \begin{array}{llll}
             F_i(X_1, \ldots, X_n) &=& X_i+H_i(X_1, \ldots, X_n), & i=1, \ldots, n,
            \end{array} 
            \label{xh}
\end{equation}            
  where $H_i \in K[X]$ has degree $D_i$ and order of vanishing $d_i$, with $ d_i\geq 2$. Let $d=\min d_i, D=\max D_i$.
  
In \cite{ABCH} an inversion algorithm for a given polynomial automorphism $F$ of the form (\ref{xh}) was proposed. Some further improvements and an estimation of the computational complexity for the algorithm 
mentioned above were discussed in \cite{AB}. 
We define a mapping $\Delta_F$ on $K[X]^n$ by $\Delta_F(P)=P\circ F-P$.
We consider the sequence $P_k$ of
polynomial maps in $K[X]^n$ defined by $P_{k}=\Delta^k_F(Id)$.
Polynomial map $F:K^n \rightarrow K^n$ is called \emph{Pascal finite} if there exists $m$ such that $P_m=0$. Pascal finite maps are invertible and the inverse map $G$ of $F$ is given by
 \begin{equation}\label{inveq}
     G(X)= \sum_{l=0}^{m-1} (-1)^l P_l(X)
 \end{equation}
 (see \cite{ABCH}, corollary 2.1). More about properties of the class of Pascal finite maps can be found in \cite{ABCH2}.
Main result of \cite{ABCH} (see theorem 3.1)  formulates an equivalent condition to
invertibility of a polynomial map and allows to  find an inverse of a given polynomial automorphism even if it is not a Pascal finite one. 

One can ask if the algorithm defined in \cite{ABCH} can be applied to mappings other than polynomial ones.
The answer is yes. It can be applied to $F=(F_1, \ldots, F_n)$, where
$F_i \in K[[X]]$ for each $i \in \{1,2, \ldots, n\}$. Here $K[[X]]$ denotes the ring of formal power series in the variables $X=(X_1, \ldots, X_n)$ with coefficients in an arbitrary field $K$.
Let us recall the following  well known result.\\

\begin{theorem}[Formal Inverse Function Theorem]
 Let $R$ be an arbitrary commutative ring, let $X=(X_1, \ldots, X_n)$. Let $F=(F_1, \ldots, F_n)$, such that $F_i \in R[[X]]$ with $F(0)=0$ and $\det{J_F}(0) \in R^*$.
 Then there exists $G=(G_1, \ldots, G_n)$ with $G_i \in R[[X]]$, such that $G(0)=0$ and $G$ is a left inverse of $F$, i.e.
 $G \circ F (X)=X$. Moreover such $G$ is uniquely determined and is also a right inverse of $F$, i.e. $F \circ G (X)=X$.
 \label{fi}
\end{theorem}

For the proof see for example \cite{E}, chapter 1. Theorem \ref{fi} states the fact that the inverse map exists. It is natural to ask about formulas for the inverse. There are many results containing answer to this question, for example Abhyankar's inversion formula (see \cite{BCW})  and others (for instance \cite{E} Theorem 3.1.1). 
Let us also mention
methods proposed for finding power series composition  by Brent and Kung (see \cite{BK1}) and Bernstein's method (see \cite{B}) and an interesting paper \cite{BDG} using  probabilistic approach to the Jacobian Conjecture.
However only some of these can be used in case of fields of positive characteristic. 

Automorphisms of formal power series rings play a significant role in algebraic geometry and as such are the subject of wide research.

The aim of this paper is to present an approach coming from the algorithmic treatment described above.
Proposed approach uses only substitution and subtraction, as a result it works for fields of arbitrary characteristic.
Power series in several variables occurring in applications often have many of their coefficients equal to zero. Our approach does not depend on any additional assumptions of this type. We establish several properties of the algorithm when applied to mappings with  polynomial inverse.

\section{Concepts needed in further considerations}\label{sec2}

Before describing our approach to inverting formal power series let us recall some notions and facts needed in futher considerations.
 
 One can define a topology on a given ring $R$. For detailed information see \cite{ZS}, chapter VIII. 2.
Let $I$ be an ideal in $R$, and consider a sequence of powers if ideal $(I^r)_{r \geq 0}$. We define \textit{I-adic topology}
of $R$ as the topology in which $I^r$ forms a basis of neighborhoods of the zero of $R$.

In considered case $R=K[[X]]$ and $I$ is the ideal generated by $X=(X_1, \ldots,X_n)$.
The ring $K[[X]]$ is a completion of the ring of polynomials $K[X]$ with respect to the $I$-adic topology. 
For every $\Phi \in K[[X]]$ we can write its decomposition in homogeneous components $\Phi = \Phi_0+\Phi_1+\ldots +\Phi_q+ \ldots $, where $\Phi_k \in I^k$ for $k \in \mathbb{N}$. The I-adic topology is exactly the topology induced by the metric $d$ given by
\[d(\Phi,\Psi)= \left\{ \begin{array}{ll}
                      C^{- ord(\Phi-\Psi)} & \Phi \neq \Psi\\
                      0 & \Phi=\Psi
                     \end{array} \right.,\] 
  for every $\Phi,\Psi \in K[[X]] $. 
  Here $C \in\mathbb{R}, C>1$ and $ord : K[[X]] \rightarrow \mathbb{Z} \cup \{ \infty \}$ is an \emph{order function}, i.e. a discrete valuation given by
 $ord(0)=\infty$ and
   $ord(\Phi)=\min \{q \in \mathbb{Z}_{ \geq 0}: \, \Phi_q \neq 0  \}$ for $\Phi \neq 0$, for every $\Phi=\Phi_0+\Phi_1+\ldots +\Phi_q+ \ldots \in K[[X]]$.

 \section{The algorithm}\label{sec3}

 Consider now $F=(F_1, \ldots, F_n) \in K[[X]]^n$.  Assume moreover that $F(0)=0$ and $F_i(X)=X_i+H_i(X)$, where $X=(X_1, \ldots, X_n)$ and the order of vanishing $ord(H_i)>1$. 
 We perform an algorithm componentwise, for each $F_i$
separately. For every $i=1, \ldots, n$ we obtain the following.

\[
\begin{array}{l}
P^i_0(X)=X_i \\
P_1^i(X)=P_0^i(F)-P_0^i(X)=H_i(X)\\
P_2^i(X)=P_1^i(F)-P_1^i(X)=H_i(F)-H_i(X) \\
\ldots \\
 P_{k+1}^i(X)=P_k^i(F)-P_k^i(X)=(P_k^i \circ F-P^i_k)(X) 
\end{array} \]

One can check that for any positive integer $m$, we have

\begin{equation}
 X= \sum_{k=0}^{m-1} (-1)^k P_k(F)+(-1)^m P_m(X).
 \label{formula}
\end{equation}

We estimate the order of vanishing in each step of the algorithm.

\begin{lemma}
  Let $F$ and $P^i_k$ be as above. Let $t_i= ord (H_i)$, for $i=1, \ldots, n$ and $t=\min t_i$.
Then $ord P^i_k \geq (k-1)(t-1)+t_i$.
\label{lowdeglem}
\end{lemma}

\begin{proof} 
 Observe that if $M=aX_1^{l_1}\ldots X_n^{l_n}$, where $a \in K$ is a monomial of degree $r=l_1+\ldots+l_n$, 
then $ord (M \circ F - M) \geq r-1+t$. Indeed, let $H_{i_0} \in I^t$. Compute $M \circ F - M$ and observe that its monomial of the smallest degree is the one of the form $a \cdot \Pi_{i \neq i_0}X_i ^{l_i} \cdot X_i^{l_{i_0}-1}\cdot H_{i_0}$. The degree of this monomial is exactly $(r-l_{i_0})+(l_{i_0}-1)+t=r-1+t$.

Every $Q\in K[[X]]$ can be written as $Q(X)=\sum_{j=s_0}^{\infty} M_{j} (X) $, for some $s_0 \in \mathbb{N}$, where 
for every $s \geq s_0$ we have $M_{s} \in I^s$ and $Q(X)-\sum_{j=s_0}^s M_{j} (X) \in I^{s+1}$. 
The order of vanishing of $Q$  is equal to the order of vanishing of $M_{s_0} \neq 0$. Consequently, if $ord(Q)=s$, then $ord(Q \circ F - Q)\geq t+s-1$.
 
We prove the estimation using induction on $k$. The thesis holds for $k=1$. Assume that $P^i_{k-1}$ has an order of vanishing $\geq (k-2)(t-1)+t_i$.
Since
$P_k^i(X) = P_{k-1}^i(F)-P_{k-1}^i(X)$, the order of vanishing of $P_k^i(X)$ is greater than or equal to
\[t+(k-2)(t-1)+t_i-1 = (k-1)(t-1)+t_i.  \]
\end{proof}

Moreover the following holds.

\begin{corollary}
 Let $K$ be a field and let $F=(F_1, \ldots, F_n) \in K[[X]]^n$ such that $F(0)=0$ and $F(X)=X+H(X)$, where the order of vanishing $t_i$ of each $H_i$ is strictly greater than one. Let $t= \min t_i$.
 Let $\{P_k\}$ be the sequence of formal power series constructed by the algorithm for $F$ and let $D \in \mathbb{Z}, \, D>1$.
For each $k=1,2,\ldots$ let $\widetilde{P_k}^D=(\widetilde{P_k^1}^D, \ldots, \widetilde{P_k^n}^D)$ denotes the polynomial mapping obtained from $P_k$ which is exactly
the sum of homogeneous summands of $P_k$ of degree $\leq D$. Let $\mu_i = \lfloor \frac{D-t_i}{t-1}+1 \rfloor +1$.
Then for every positive integer $q$ we have

\[ \sum_{k=0}^{\mu_i-1}(-1)^k \widetilde{P^i_k}^D (X_1, \ldots, X_n) = \sum_{k=0}^{\mu_i-1+q}(-1)^k \widetilde{P^i_k}^D (X_1, \ldots, X_n).\]
 \label{stepnumberlem}
\end{corollary}
\begin{proof} 
 According to lemma \ref{lowdeglem} the order of vanishing of $P^i_k$ is equal or greater than
 $(k-1)(t-1)+t_i$. Hence if $k>\frac{D-t_i}{t-1}+1 $, the order of vanishing of $P_k$ is greater than $D$.
 So $\mu_i$ is the smallest such $k$ and we get the thesis.
\end{proof}

  \begin{remark}
  Let $K$ be a field and let $F=(F_1, \ldots, F_n) \in K[[X]]^n$ ba as above. Let $(P_n)$ be the sequence of formal power series constructed by the algorithm for $F$. Due to corollary \ref{stepnumberlem}, the sequence $(A_m)$ of elements in
  $K[[X]]^n$ of the form $A_m = \sum_{k=0}^{m-1} (-1)^kP_k(X)$ converges to some $A \in K[[X]]^n$.
 \end{remark}

 \begin{theorem}
  Let $K$ be a field and let $F=(F_1, \ldots, F_n) \in K[[X]]^n$ such that $F(0)=0$ and $F(X)=X+H(X)$, where the order of vanishing $t_i$ of each $H_i$ is strictly greater than one. Let $t=\min t_i$ and $G:=F^{-1}$
be the formal inverse of $F$.
 Let $\{P_k\}$ be the sequence of formal power series constructed by the algorithm for $F$.
 Denote $A_m = \sum_{k=0}^{m-1} (-1)^kP_k(X)$ and $A= \sum_{k=0}^{\infty} (-1)^kP_k(X)$. Then $A=G$.
  \label{fithm}
 \end{theorem}
\begin{proof} 
 We need to show that for every $i=1, \ldots, n$ we have $\sum_{k=0}^{\infty} (-1)^kP_k^i(X)=G_i(X)$ or equivalenty (since $G_i(F)=X_i$)
 that $\sum_{k=0}^{\infty} (-1)^kP_k^i(F)=X_i$. Due to (\ref{formula}) for any positive integer $m$,

\[ \sum_{k=0}^{m-1} (-1)^k P_k^i(F_1,\dots,F_n)= X_i -(-1)^m P_m^i(X_1,\dots,X_n) = X_i +(-1)^{m-1} P_m^i(X_1,\dots,X_n). \]
So it is enough to prove that $(-1)^{m-1} P_m^i(X_1,\dots,X_n)$ tends to zero.
And this is exactly the case. By lemma \ref{lowdeglem} the order of vanishing of $P^i_m$ is greater than or equal to $(m-1)(t-1)+t_i$.
Since $t_i>1$ for each $i=1, \ldots, n$, then $t-1>0$ and the order of vanishing of $P^i_m$ tends to infinity when $m \rightarrow \infty$.
\end{proof}

\section{Mappings with polynomial inverse}\label{sec4}

Let us now consider a special case. Assume that we are dealing with $F \in K[[X]]^n$ such that
its formal inverse $G:=F^{-1}$, where $G=(G_1,\ldots,G_n)$ is a \textit{polynomial mapping} of degree $D$. Then we can recover a certain version 
of theorem 3.1 in \cite{ABCH}.

\begin{corollary} \label{newsymthm}
 Let $K$ be a field and let $F=(F_1, \ldots, F_n) \in K[[X]]^n$ such that $F(0)=0$ and $F(X)=X+H(X)$, where the order of vanishing $t_i$ of each $H_i$ is strictly greater than one. Let $t= \min t_i$. Assume that  
the formal inverse $G=F^{-1}$ is a polynomial mapping of degree $D$.
 Let $\{P_k\}$ be the sequence of formal power series constructed by the algorithm for $F$.
 Then the following holds.
 
 \begin{enumerate}
  \item[1)] For $i=1,\ldots, n$ and  $\mu_i = \lfloor \frac{D-t_i}{t-1}+1 \rfloor +1$, we have
   \begin{equation}\label{sc}
    \forall m \in \mathbb{Z} \qquad m \geq \mu_i \quad \Rightarrow \quad \sum_{k=0}^{m-1}(-1)^k P^i_k(X)=G_i(X)+R^i_m(X), 
   \end{equation}
   where $G_i(X)$ is a polynomial of degree $ \leq D$ independent of $m$,
  and $R^i_m(X) \in K[[X]]$ satisfies $R^i_m(F)=(-1)^{m+1}P^i_m(X)$.
  \item[2)]  Moreover the inverse $G$ of $F$ is given by

  \[G_i(Y_1, \ldots, Y_n)=\sum_{k=0}^{m-1}(-1)^k\widetilde{P^i_k}^D(Y_1, \ldots, Y_n), \, i=1, \ldots, n,\]
  where $\widetilde{P^i_k}^D$ is the sum of homogeneous summands of $P^i_k$ of degree $ \leq D$ and $m$ is any integer grater than $ \frac{D-t_i}{t-1}+1$.
 \end{enumerate}

\end{corollary}

\begin{proof} 
 1) By applying (\ref{formula}), we obtain, for any positive integer $m$

\begin{equation*}
X_i= \sum_{k=0}^{m-1} (-1)^k P_k(F_1,\dots,F_n)+(-1)^m P_m(X_1,\dots,X_n).
\end{equation*}

\noindent Since $X_i=G_i(F_1,\dots,F_n)$, we obtain the
following equality of polynomials in the variables $Y_1,\dots,Y_n$.

\begin{equation*}
\sum_{k=0}^{m-1} (-1)^k P_k(Y_1,\dots,Y_n)=  G_i(Y_1,\dots,Y_n)  -(-1)^m P_m(G_1(Y_1,\dots,Y_n),\dots,G_n(Y_1,\dots,Y_n)),
\end{equation*}

\noindent Denote

$$R_m^i(Y_1,\dots,Y_n):=-(-1)^m P_m^i(G_1(Y_1,\dots,Y_n),\dots,G_n(Y_1,\dots,Y_n)).$$

\noindent We obtain $\sum_{k=0}^{m-1}(-1)^k P^i_k(Y)=G_i(Y)+R^i_m(Y)$, where $Y=(Y_1, \ldots, Y_n)$. Now, $G_i$ is a polynomial of degree at most $D$ in
$Y_1, \dots, Y_n$. 
According to lemma \ref{lowdeglem} for an integer $m$ such that $m > \frac{D-t_i}{t-1}+1$, $P_m^i$ is a
formal power series in the variables $X_1,\dots,X_n$  of order of vanishing bigger than $D$, hence the order of vanishing of $
P_m(G_1(Y_1,\dots,Y_n),\dots,G_n(Y_1,\dots,Y_n))$ in the variables $
Y_1,\dots,Y_n$ is bigger than $D$. Therefore, the sum of homogeneous
summands of degrees not bigger than $D$ in the right-hand side of the
equality above is exactly $\sum_{k=0}^{m-1} (-1)^k \widetilde{P
_k^i}^D(Y_1,Y_2,\dots,Y_n)$.

2) Again due to (\ref{formula}) we obtain, for any positive integer $m$,

\[ 
\begin{array}{rl}
 X_i=& \sum_{k=0}^{m-1} (-1)^k P_k^i(F)+(-1)^m P_m^i(X) \\
=&G_i(F)+R_m^i(F)+(-1)^m P_m^i(X) \\
=& G_i(F)+(-1)^{m+1} P_{m}^i(X)+(-1)^m P_m^i(X)\\
=& G_i(F). 

\end{array}\]
Hence the formal  inverse of $F$ is $G$.
\end{proof}

\noindent From now on, we will refer to condition (\ref{sc}) in the theorem above as the \emph{symmetry condition}. \\

\section{Examples}\label{sec5}

\begin{example} Consider $G: \mathbb{R} \rightarrow \mathbb{R}$,  $G(Y)=Y+Y^2$. We have $G(0)=0$, $G'(Y)=1+2Y$ and $G'(0)=1$. The inverse function is a power series of the form
\[Y=F(X)=\sum_{j=0}^{\infty} (-1)^j \binom{2j}{j} \frac{X^{j+1}}{j+1} = X-X^2+2X^3-5X^4 +14 X^5 - 42 X^6   +O(X^7).\]  Here $n=1$,$D=t=2$ and $\mu =2$.
We perform the algorithm for $F$ and obtain the following.
\[
\begin{array}{l}
P_0(X)=X,\\
P_1(X)= -X^2 + 2 X^3 - 5 X^4 + 14 X^5 - 42 X^6  +O(X^7)\\
 P_2(X) = 2 X^3 - 11 X^4 + 52 X^5 - 238 X^6 + O(X^7)
\end{array} \]
Observe that the order of vanishing of each $P_k$ is $k+1=(k-1)(d-1)+d$.
Moreover since $F$ has a polynomial inverse, the symmetry condition holds. Indeed we have the following.

\[
\begin{array}{rl}
P_0(X)-P_1(X)=&X+X^2 - 2 X^3 + 5 X^4 - 14 X^5 + 42 X^6 + O(X^7) \\
G(X) =& \widetilde{P_0}^2(X)-\widetilde{P_1}^2(X)=X+X^2\\
R_2(X) =& - 2 X^3 + 5 X^4 - 14 X^5 + 42 X^6 + O(X^7)\\
R_2(Y)=&(-1)^{2+1}P_2\circ G(Y)=-P_2(Y+Y^2)
\end{array} \]

Consider now $G: \mathbb{R} \rightarrow \mathbb{R}$,  $G(Y)=Y+Y^2$, which is the inverse of the mapping $F$ given  above.
We perform the algorithm for $G$ and obtain the following.
\[
\begin{array}{rl}
Q_0(Y)=&Y\\
Q_1(Y)= &Y^2\\
 Q_2(Y) = &Y^4+2Y^3 \\
 Q_3 (Y)=  &Y^8+4Y^7+8Y^6+10Y^5+6Y^4\\
 Q_4(Y) =  & Y^{16}+8Y^{15}+32Y^{14}+84Y^{13}+162Y^{12}+244Y^{11}\\
 &+298Y^{10}+302Y^9+254Y^8+172Y^7+86Y^6+24Y^5
\end{array} \]

\label{sp}
 \end{example}

 \begin{remark} Considering power series with a polynomial inverse is interesting since it raises the natural question of whether there is a counterexample to the Jacobian Conjecture.
   Suppose that we have a map $F \in K[[X]]^n$ for which there exists a formal inverse.  We will use the notation established in the corollary \ref{newsymthm}. 
Let $G \in K[Y]^n$ be its inverse, which is a polynomial map of degree $D$. We have $G(Y)=Y+ L(Y)$. Denote $d_i=ord (L_i)$, $d=\min d_i$. Of course $G$ is not Pascal finite. If it was, $F$ would be too (see \cite{ABCH2}, theorem 3.1). 
By corollary \ref{newsymthm} the symmetry condition holds for $F$.
The proof of point 1) depends on the degree of the inverse. 
As a result, we cannot expect the symmetry condition to hold for $G$ (see \cite{ABCH}, theorem 3.1).

 \end{remark}

\begin{example}
 Consider $F: \mathbb{R} \rightarrow \mathbb{R}$, where $F(X)=\sin{X}$. We have
 \[  \begin{array}{l}
   F(X)= \sum_{j=0}^{\infty} \frac{(-1)^j}{(2j+1)!}X^{2j+1} = X - \frac{1}{6}X^3 + \frac{1}{120}X^5 - \frac{1}{5040}X^7 + \frac{1}{362880}X^9 - \frac{1}{39916800}X^{11} \\
   \\
    + \frac{1}{6227020800}X^{13} - \frac{1}{1307674368000}X^{15} + \frac{1}{355687428096000}X^{17} - \frac{1}{121645100408832000}X^{19} +\\
    \\
 +\frac{1}{51090942171709440000}X^{21} - \frac{1}{25852016738884976640000}X^{23} + \frac{1}{15511210043330985984000000}X^{25} \\
 \\
 - \frac{1}{10888869450418352160768000000}X^{27} + \frac{1}{8841761993739701954543616000000}X^{29} + O(X^{30})
 \end{array}
 \]
 
 We execute the algorithm and obtain partial sums $A_m$ as following.
 \end{example}

 \begin{table}[!t]
 \caption{Partial sums $A_m$ obtained by the algorithm.}\label{tab:stabil}
{\begin{tabular}{@{}llll@{}}
\hline
$m$ & $A_m$ \\
\hline
1&$X$
     \\ \hline
     \\
     2& $X + \frac{1}{6}X^3 - \frac{1}{120}X^5 + \frac{1}{5040}X^7 - \frac{1}{362880}X^9 + \frac{1}{39916800}X^{11} - \frac{1}{6227020800}X^{13} +\frac{1}{1307674368000}X^{15} $\\
     \\ 
     &$- \frac{1}{355687428096000}X^{17} +
     \frac{1}{121645100408832000}X^{19} - \frac{1}{51090942171709440000}X^{21}$\\
     \\ 
     &$+\frac{1}{25852016738884976640000}X^{23}- \frac{1}{15511210043330985984000000}X^{25}  +\frac{1}{10888869450418352160768000000}X^{27}$\\
     \\ 
 &$- \frac{1}{8841761993739701954543616000000}X^{29} + O(X^{30})$\\
 \\ \hline
 \\
 3& $X + \frac{1}{6}X^3 + \frac{3}{40}X^5 - \frac{25}{1008}X^7 + \frac{89}{17280}X^9 - \frac{5371}{5702400}X^{11} + \frac{90071}{566092800}X^{13} 
  - \frac{2487217}{100590336000}X^{15}$\\
  \\
  &$+ \frac{181808603}{50812489728000}X^{17}- \frac{3501606701}{7155594141696000}X^{19} + \frac{171190842799}{2688996956405760000}X^{21}$\\
  \\
   &$- \frac{68118191390719}{8617338912961658880000}X^{23}+ \frac{42303813823301}{44960029111104307200000}X^{25}- \frac{167768077105452763}{1555552778631193165824000000}X^{27} $\\
   \\
  &$+\frac{105595092426069327869}{8841761993739701954543616000000}X^{29} + O(X^{30})$\\
 & \\ \hline
 \\
 4& $X + \frac{1}{6}X^3 + \frac{3}{40}X^5 + \frac{5}{112}X^7 - \frac{175}{3456}X^9 + \frac{11057}{380160}X^{11} - \frac{24559}{1797120}X^{13} + \frac{849089}{145152000}X^{15}$\\
 \\
 &$- \frac{64936446307}{27360571392000}X^{17} +\frac{2215609252019}{2385198047232000}X^{19} - \frac{18638240496367}{52725430517760000}X^{21} + \frac{16240703065022129}{123693859994664960000}X^{23}$\\
 \\
 &$- \frac{238326243945432103}{4995558790122700800000}X^{25} +\frac{41248891834131025189}{2427841571999632588800000}X^{27} 
  - \frac{3089627793964198212305399}{520103646690570703208448000000}X^{29}$\\
  \\
  &$+ O(X^{30})$\\
 & \\ \hline
 \\
 5& $X + \frac{1}{6}X^3 + \frac{3}{40}X^5 + \frac{5}{112}X^7 + \frac{35}{1152}X^9 - \frac{2513}{25344}X^{11} + \frac{629467}{5391360}X^{13} - \frac{430757}{4147200}X^{15}$\\
 \\
 &$+ \frac{1918781351}{23688806400}X^{17} - \frac{978413543059}{16679706624000}X^{19} + \frac{101839919706341}{2510734786560000}X^{21} - \frac{197803541540754287}{7276109411450880000}X^{23}$\\
 \\
 &$+ \frac{2299085857670705567}{129071853907476480000}X^{25}- \frac{27844882444217079689291}{2427841571999632588800000}X^{27}+ \frac{318927153969304780440241}{43809052365860036935680000}X^{29}$\\
 \\
 &$+ O(X^{30})$\\
 & \\ \hline
\end{tabular}}
\end{table}

\newpage
  And finally we obtain the inverse.
  \[ \begin{array}{l}
     G(X) = X + \frac{1}{6}X^3 + \frac{3}{40}X^5 + \frac{5}{112}X^7 + \frac{35}{1152}X^9 + \frac{63}{2816}X^{11} + \frac{231}{13312}X^{13} + \frac{143}{10240}X^{15} + \frac{6435}{557056}X^{17} \\
     \\
  \frac{12155}{1245184}X^{19} + \frac{46189}{5505024}X^{21} + \frac{88179}{12058624}X^{23} + \frac{676039}{104857600}X^{25} + \frac{1300075}{226492416}X^{27} + \frac{5014575}{973078528}X^{29} + O(X^{30}) 
    \end{array} \]

  \begin{remark}
 
  According to theorem \ref{fithm} an inverse $G$ of a mapping $F(X)=X+H(X)$ is a limit $A$ of the sequence of partial sums $A_m=(A_m^1, \ldots, A_m^n)$, where $A_m^i = \sum_{k=0}^{m-1} (-1)^kP_k^i(X)$.
 Let $t_i = ord(H_i)$ denote the order of vanishing of  $H_i$, for $i=1, \ldots, n$ and $t=\min t_i$.
By lemma \ref{lowdeglem} the formal power series $P^i_k$ has  order of vanishing $\geq (k-1)(t-1)+t_i$.
That is why in each sequence  of power series $(A_m^i)_{m \geq 0}$ for $i=1, \ldots, n$ one can observe stabilization of terms, starting from those of small degree (see for example table \ref{tab:stabil}). Precisely, since $ord(P^i_{m-1}) \geq (m-2)(t-1)+t_i$ then all terms of degree strictly smaller than $(m-2)(t-1)+t_i$ appearing in $A^i_m$ are already fixed (i.e. the same as in $G=A$).
  \label{stabilization}
  \end{remark}
  
  \begin{remark}
   When performing calculations on power series with any Computer Algebra System one has to set some precision. Due to remark \ref{stabilization} one can compute number of steps of algorithm that need to be performed in order to obtain an inverse mapping $G$ with a given precision. Precisely, in order to obtain $G=(G_1,\ldots, G_n)$ with precision $O(x_1,\ldots, x_n)^s$, one needs to perform as many steps as it takes to observe stabilization of terms of degree smaller or equal to $s$ in each $A^1_m, \ldots, A^n_m$. This happens precisely when $m$ is such that
   \[\min\{ord(P^1_{m-1}), \ldots, ord(P^n_{m-1})\} = (m-2)(t-1)+t > s .\]
   Hence $m>2+\frac{s-t}{t-1}$ and $m_0=\lfloor 2+\frac{s-t}{t-1}\rfloor +1$ is minimal such $m$.
   \label{precision}
   \end{remark}

  \begin{example} \label{ex:sin:cos} Consider
  \begin{equation*}
   F=(F_1,F_2): \mathbb{R}^2 \rightarrow \mathbb{R}^2, \, F(X,Y)=\Big( \sin(X+Y) - Y ,\cos(XY) - 1  + Y \Big).
  \end{equation*}
 Here 
 \[ 
  \begin{array}{l}
    F_1(X,Y) = X - \frac{1}{6}X^3 - \frac{1}{2}X^2Y - \frac{1}{2}XY^2 - \frac{1}{6}Y^3 + \frac{1}{120}X^5 + \frac{1}{24}X^4Y + \frac{1}{12}X^3Y^2 + \frac{1}{12}X^2Y^3 \\
    \\
    + \frac{1}{24}XY^4  + \frac{1}{120}Y^5  - \frac{1}{5040}X^7 
 - \frac{1}{720}X^6Y - \frac{1}{240}X^5Y^2 - \frac{1}{144}X^4Y^3 - \frac{1}{144}X^3Y^4 - \frac{1}{240}X^2Y^5\\
 \\
 - \frac{1}{720}XY^6 - \frac{1}{5040}Y^7  + \frac{1}{362880}X^9 + \frac{1}{40320}X^8Y + \frac{1}{10080}X^7Y^2 + \frac{1}{4320}X^6Y^3 + \frac{1}{2880}X^5Y^4 \\
 \\
 + \frac{1}{2880}X^4Y^5 + \frac{1}{4320}X^3Y^6 + \frac{1}{10080}X^2Y^7  
+ \frac{1}{40320}XY^8 + \frac{1}{362880}Y^9 + O(X, Y)^{10}
  \end{array}
\]
 and
 \[\begin{array}{l}
    F_2(X,Y)=Y - \frac{1}{2}X^2Y^2 + \frac{1}{24}X^4Y^4 + O(X, Y)^{10}.
   \end{array}\]

 Let us assume that we want to find an inverse $G(X,Y)=(G_1(X,Y), G_2(X,Y))$ with precision $O(X, Y)^{10}$. Using the notation established in theorem \ref{fithm} and remark \ref{precision} we have $t_1=3, t_2=4$ and $t=3, s=10$. Hence $m > \frac{11}{2}$ and $m_0=6$. Consequently we perform six steps of the algorithm giving $P_0, \ldots, P_5$  in order to compute $A_6$ and obtain
 \[
   \begin{array}{l}
G_1(X,Y) = X + \frac{1}{6}X^3 + \frac{1}{2}X^2Y + \frac{1}{2}XY^2 + \frac{1}{6}Y^3 + \frac{3}{40}X^5 + \frac{3}{8}X^4Y + \frac{3}{4}X^3Y^2 + \frac{3}{4}X^2Y^3 \\ 
\\
 + \frac{3}{8}XY^4 + \frac{3}{40}Y^5 + \frac{1}{4}X^4Y^2  + \frac{1}{2}X^3Y^3 + \frac{1}{4}X^2Y^4 + \frac{5}{112}X^7 + \frac{5}{16}X^6Y + \frac{15}{16}X^5Y^2 + \frac{25}{16}X^4Y^3 \\
 \\
 + \frac{25}{16}X^3Y^4 + \frac{15}{16}X^2Y^5 + \frac{5}{16}XY^6 + \frac{5}{112}Y^7 + \frac{13}{48}X^6Y^2 + \frac{7}{6}X^5Y^3 + \frac{47}{24}X^4Y^4 + \frac{19}{12}X^3Y^5\\
 \\
 + \frac{29}{48}X^2Y^6 + \frac{1}{12}XY^7  + \frac{35}{1152}X^9 + \frac{35}{128}X^8Y + \frac{35}{32}X^7Y^2  
 + \frac{269}{96}X^6Y^3 + \frac{285}{64}X^5Y^4 + \frac{269}{64}X^4Y^5 \\
 \\
 + \frac{245}{96}X^3Y^6 + \frac{35}{32}X^2Y^7 + \frac{35}{128}XY^8 + \frac{35}{1152}Y^9 + O(X, Y)^{10}
   \end{array} \]
 \[ 
 \begin{array}{l}
  G_2(X,Y) = Y + \frac{1}{2}X^2Y^2 + \frac{1}{6}X^4Y^2 + \frac{1}{2}X^3Y^3 + \frac{1}{2}X^2Y^4 + \frac{1}{6}XY^5 + \frac{1}{2}X^4Y^3 + \frac{4}{45}X^6Y^2\\
  \\
  + \frac{11}{24}X^5Y^3  + \frac{11}{12}X^4Y^4 + \frac{37}{36}X^3Y^5 + \frac{7}{12}X^2Y^6 + \frac{19}{120}XY^7 + \frac{1}{72}Y^8 + \frac{1}{3}X^6Y^3
  + \frac{5}{4}X^5Y^4\\
  \\
  + \frac{3}{2}X^4Y^5 + \frac{7}{12}X^3Y^6 + O(X, Y)^{10}
 \end{array} \]
  \end{example}

\section*{Acknowledgement}

This research was supported by the AGH University of Krakow within subsidy from Polish Ministry of Science and Higher Education (grant no. 16.16.420.054).
The author has no other relevant financial or non-financial interests to disclose.

\end{document}